\documentstyle[amssymb,amsfonts]{amsart}

 at 10 true pt

\def\beq{\begin{equation}}
\def\eeq{\end{equation}}
\def\barray{\begin{eqnarray*}}
\def\earray{\end{eqnarray*}}
\def\be#1{ \begin{equation}\label{#1} }

\def\bas{\begin{align*}}
\def\eas{\end{align*}}
\def\bi{\begin{itemize}}
\def\ei{\end{itemize}}

\def \endprf{\hfill  {\vrule height6pt width6pt depth0pt}\medskip}
\def\emph#1{{\it #1}}
\def\textbf#1{{\bf #1}}

\def\hs{\hfill $\square$}

\theoremstyle{plain}
  \newtheorem{theorem}[subsection]{Theorem}

  \newtheorem{lemma}[subsection]{Lemma}
  \newtheorem{cor}[subsection]{Corollary}

\theoremstyle{remark}
  \newtheorem{remark}[subsection]{Remark}

\theoremstyle{definition}

\include{psfig}

\begin{document}

\title{  Olson's theorem for cyclic groups}

\author{  Van Vu}
\address{Department of Mathematics, UCSD, La Jolla, CA 92093-0112}
\email{vanvu@@ucsd.edu}

\thanks{V. Vu is an A. Sloan  Fellow and is supported by an NSF Career Grant.}

\begin{abstract} Let $n$ be a large number. A subset $A$ of $Z_n$
is complete if $S_A = Z_n$, where $S_A$ is the collection of the
subset sums of $A$. Olson proved that if $n$ is prime and $|A|>
2n^{1/2} $,  then $S_A$ is complete. We show that  a similar
result for the case when $n$ is a composite number, using a
different approach.

\end{abstract}

\maketitle

\section {Introduction}

Let $G$ be an additive  group. For a subset $A \subset G$, we
denote by $S_A$ the collection of the subset sums of $A$

$$S_A =\{ \sum_{x \in B} x | B \subset A, |B| < \infty \}. $$

Following \cite{erd}, we  say that $A$ is {\it complete} (with
respect to $G$) if $S_A=G$; in other words, every element of $G$
can be represented as a sum of different elements of $A$.

In this short note, we investigate the case when $G=Z_n$, the
cyclic group of order $n$,  where $n$ is a large positive integer.

A well-known result of Olson \cite{Olson}, answering a question of
Erd\"os and Heilbronn, shows that if $n$ is a prime and $|A| >
2n^{1/2}$, then $A$ is complete.

\begin{theorem} \label{theo:olson}  If $n$ is a prime and $A$ is a
subset of $Z_n$ with cardinality larger than $2n^{1/2}$, then $A$
is complete.
\end{theorem}

The bound is sharp. To see this observe that if the sum of the
elements in $A$ (viewed as integers between $1$ and $n-1$) is less
than $n$, then $A$ is not complete.

We  extend this result for the case when $n$ is a composite
number. Our result is

\begin{theorem} \label{theo:main}  There is a constant $C$ such that the following holds.
Let $n$ be a  sufficiently large positive integer  and  $A$ be a
subset of  $Z_n$, where $|A| \ge Cn^{1/2}$ and the elements of $A$
are co-primes with $n$. Then  $A$ is complete.
\end{theorem}

\begin{remark} \label{remark:diderich} The assumption that the
elements of $A$ are co-primes with $n$ is necessary. For instance,
if $n$ is divisible by 3 then it is possible to have  an
incomplete set of size $n/3$. Without the co-prime assumption, the
problem of bounding $|A|$ is known as Diderrich's problem. It has
been proved that the sharp bound for $|A|$ is $p+n/p-2$, where $p$
is the smallest prime divisor of $n$ (see \cite{Lipkin} for the
case of cyclic groups and \cite {Gao} for the general case of
arbitrary abelian groups).
\end{remark}

In the current proof, the constant $C$ in Theorem \ref{theo:main}
is fairly large. However, we believe that the constant $C$ in
Theorem \ref{theo:main} can be set to (the asymptotically optimal
value) $2+o(1)$.

It seems to be of interest  to investigate the general case. Given
a finite abelian group $G$, one would like to find a parameter
$f(G)$  so that if $A$ is a set of at least $f(G)$ primitive
elements, then $A$ is complete. This problem can be seen as a
variant of Diderrich's problem. On the other hand, by comparing
the bounds in Theorem \ref{theo:main} and Remark
\ref{remark:diderich}, it is plausible that the answer would be
quite different. In fact, we think  that the nature of this
problem is closer  to  that of Olson's than to Diderrich's.

\vskip2mm

{\bf Notation.}  In the whole paper, we understand that the
elements of a set are different. If there are possible repetitions
we use the phrase multi-set instead.

\section{Lemmas}

\noindent For a set $A$ of integers and a positive integer $l \le
|A|$, let $l^{\ast}A$ denote the set of sums of $l$ different
elements of $A$

$$l^{\ast}A = \{a_1+ \dots a_l| a_i \in A, a_i \neq a_j\}. $$

\noindent  Denote by $[n]$ the set $\{1,2, \dots, n\}$. In
\cite{szemvu1}, Szemer\'edi and the author proved the following
theorem.

\begin{theorem} \label {theo:szemvu1} There are positive constants $C$
and $c$ such that the following holds. Let $l$, $n$ be positive
integers and  $A$ be a subset of $[n]$ such that $|A|/2 \ge l$ and
$l|A| \ge Cn$. Then $l^{\ast}A$ contains an arithmetic progression
of length $cl|A|$.
\end{theorem}

\noindent Since $l^{\ast}A \subset S_A$, this theorem implies the
following corollary.

\begin{cor} \label {cor:szemvu1} There is a positive constant $C$
such that the following holds. For every sufficiently large
integer $n$ and a subset $A$ of $[n]$ of cardinality at least $C
n^{1/2}$, $S_A$ contains an arithmetic progression of length $n$.
\end{cor}

\begin{remark} \label{remark:FS} The bounds on both $|A|$ and the length of the arithmetic progression is
sharp, up to constant factors. Freiman \cite{Fre} and
 S\'ark\"ozy \cite{Sar}, independently,  showed that  the same statement holds under the stronger assumption that
 $|A| \ge C \sqrt {n \log n}$.   \end{remark}

\noindent We also need the following simple lemma.

\begin{lemma} \label{lemma:lemma2} Let $n$ be a positive integer and  $A$ be a multi-set of
$n$  integers co-prime to $n$. Then $S_A$ contains every residue
modulo $n$.
\end{lemma}

{\bf \noindent Proof of Lemma \ref{lemma:lemma2}.} Assume that
$a_1, a_2, \dots , a_n$ are the elements of $A$. We are going to
prove, by induction, that $|S_{A_i}| \ge i$, where $A_i=\{a_1,
\dots, a_i \}$. The case $i=1$ is trivial. Assume that the
statement holds for $i-1$. Let $b_1, \dots, b_{i-1}$ be $i-1$
different elements (modulo $n$) of $S_{A_{i-1}}$. Since  the
statement is invariant under dilation,
 we can assume that $a_i=1$. Consider the elements

 $$b_1, \dots, b_{i-1}, 1, 1+ b_1, \dots, 1+ b_{i-1}. $$

 At least $i$ of the above must be different (modulo $n$) and this concludes
 the proof. \hs

 \section {Proof of Theorem \ref{theo:main}.} Assume that $A$ has
 at least
 $2 \lceil C  n^{1/2} \rceil $ elements, where $C$ is the
 constant in Corollary \ref{cor:szemvu1}. For convenience, we  think of the elements of $A$
 as  positive integers between one and $n-1$. We are going to
 prove that $S_A$  contains every residue modulo $n$.

 Let $A'$
be a subset of $A$ of $ \lceil Cn^{1/2} \rceil $ elements. Apply
Corollary \ref{cor:szemvu1} to $A'$ to get an arithmetic
progression $P'$ of length $n$. If the difference $d'$ of $P'$ is
co-prime to $n$, then $P'$ contains every residue modulo $n$ and
we are done. If $d'$ is not co-prime to $n$, set $d= \gcd (d',n)$.

Since the largest element in $S_{A'}$ is less than $Cn^{3/2}$, $d$
is less than $Cn^{1/2}$. The set $B= A \backslash A'$ has at least
$ \lceil Cn^{1/2} \rceil
> d $ elements, each of which is  co-prime to $n$ (and thus co-prime to $d$). By Lemma
\ref{lemma:lemma2}, we conclude that $S_{B'}$ contains every
residue modulo $d$.

The set $S_{A'} +S_{B}$ thus contains every residue modulo $n$.
But this set is clearly a subset of $S_A$, completing the proof.
\hs


Notice that the proof requires that the elements of $B$ are
co-prime to $n$; but for $A'$, it is enough to assume that its
elements are non-zero modulo $n$.

{\it Acknowledgement.} We would like to thank T. Tao for reading
the manuscript.

\end{document}